\documentclass[english,reqno]{amsart}
\usepackage{bookman}
\usepackage[T1]{fontenc}
\usepackage[latin1]{inputenc}
\IfFileExists{url.sty}{\usepackage{url}}
                      {\newcommand{\url}{\texttt}}

\makeatletter

 \theoremstyle{plain}
 \theoremstyle{plain}    
 \newtheorem{thm}{Theorem} 
 \theoremstyle{remark}    
 \newtheorem*{claim*}{Claim}
 \theoremstyle{remark}
 \newtheorem*{rem*}{Remark}

\newcommand{\NN}{\mathbb {N}}

\newcommand{\RR}{\mathbb {R}}

\newcommand{\Us}{\mathcal {U}}

\newcommand{\St}{\mathrm{St}\,}

\allowdisplaybreaks[2]

\usepackage{babel}
\makeatother
\begin{document}

\title{On removing one point from a compact space}

\author{Gady Kozma}

\email{gady@ias.edu, gadykozma@hotmail.com}

\begin{abstract}
If $B$ is a compact space and $B\setminus\{ pt\}$ is Lindel\"of then
$B^{\kappa}\setminus\{\overrightarrow{pt}\}$ is star-Linedl\"of
for every $\kappa$. If $B\setminus\{ pt\}$ is compact then $B^{\kappa}\setminus\{\overrightarrow{pt}\}$
is discretely star-Lindel\"of. In particular, this gives new examples
of Tychonoff discretely star-Lindel\"of spaces with unlimited extent.
\end{abstract}
\maketitle

\section{Introduction}

A topological space $X$ is called star-Lindel\"of if for every open
cover $\Us$ there exists a countable set $F$ such that $\St(F,\Us)=X$
where\[
\St(F,\Us):=\bigcup\{ U\in\Us\,:\, U\cap F\neq\emptyset\}\]
or in other words, one may extract from $\Us$ a subcover which can
be written as a countable collection of {}``stars'', each centered
around a point of $F$. It is an interesting and widely researched
class (see \cite{M98} for a survey) containing many spaces, including
all Lindel\"of spaces, all separable spaces and all countably compact
spaces, for example $\omega_{1}$. A space is discretely star-Lindel\"of
if the set $F$ can also be taken to be closed and discrete (this
property is called {}``in discrete web'' in \cite{YY99} and \cite{M00}).
Generally, products of star-Lindel\"of spaces are not star-Lindel\"of, and
even a product of a star-Lindel\"of space with a compact space need
not be star-Lindel\"of (see \cite{BM01}). In light of the results
stated above (and even more so, in light of their proofs below), one
might be tempted to think that the following might be true:

\theoremstyle{plain} 

\newtheorem*{untrue}{Untrue}

\begin{untrue}

For every compact space $B$ such that $B\setminus\{ pt\}$ is star-Lindel\"of,
so is $B^{\kappa}\setminus\{\overrightarrow{pt}\}$.

\end{untrue}

This is not true even for the product of 2 spaces. Let $B$ be a compact
space and $0\in B$ some point such that $B\setminus\{0\}$ is a star-Lindel\"of
space with uncountable extent (see below for the definition). Then
a proof very similar to lemma 2.3 of \cite{BM01} shows that $(B\cup[0,\omega_{1}])^{2}\setminus\{(0,\omega_{1})\}$
is not star-Lindel\"of. It is also possible to identify the two points
$0$ and $\omega_{1}$ and get an example where the point removed
is $(0,0)$.

\theoremstyle{plain} 

\newtheorem{question}{Question}

\begin{question}

Is there a star-Lindel\"of-like property T such that if $B$ is compact
and $B\setminus\{ pt\}$ is T then $B^{\kappa}\setminus\{\overrightarrow{pt}\}$
is T for every $\kappa$?%
\footnote{It can be arguably claimed that this question is phrased too vaguely
to be answered in the negative...%
}

\end{question}

I must add that it is not immediate that the two results of this paper
cannot be joined. Thus I have

\begin{question}

Is it true that if $B$ is compact and $B\setminus\{ pt\}$ is Lindel\"of
then $B^{\kappa}\setminus\{\overrightarrow{pt}\}$ is discretely star-Lindel\"of
for every $\kappa$?

\end{question}

One of the motivations to discuss these point-removed-from-a-product
spaces is that they tend to have large extent. The extent of a space
$X$, denoted by $e(X)$, is the supremum of the cardinalities of
closed discrete subspaces of $X$. The connection between the star-Lindel\"ofness or
the star-Lindel\"of number%
\footnote{See definition on page \pageref{def_starlinnum}.%
} of a space and its extent is a natural one --- actually, the star-Lindel\"of
number is sometimes called the weak extent (see \cite{M02}). In \cite{B98}
this connection is discussed but only partial results are obtained.
Much stronger results were obtained in \cite{M02} where Tychonoff
(i.e.~$T_{3\frac{1}{2}}$) star-Lindel\"of spaces of unbounded extent
were constructed, and it was shown that the extent of a $T_{4}$ star-Lindel\"of
space is $\leq2^{\aleph_{0}}$ (clearly, the extent of a metric star-Lindel\"of
space is at most countable). In \cite{M00} it was further shown that
the existence of a $T_{4}$ star-Lindel\"of space with uncountable
extent is consistent with ZFC. The following are still open, though:

\begin{question}

Is the existence of a $T_{4}$ star-Lindel\"of with a closed discrete
subset of cardinality $2^{\aleph_{0}}$ consistent?

\end{question}

I was informed by R. Levy that something quite close may be demonstrated:
a model in which $2^{\aleph_{0}}$ is a limit cardinal and there exists
a $T_{4}$ star-Lindel\"of space $X$ (actually, a separable space)
with closed discrete subsets $D^{\lambda}$ satisfying $|D^{\lambda}|\rightarrow2^{\aleph_{0}}$,
so that $e(X)=2^{\aleph_{0}}$.

\begin{question}

Is the existence of a $T_{4}$ star-Lindel\"of space with uncountable
extent indeed independent from ZFC?

\end{question}

The same connection between extent and separation can be questioned
for discretely star-Lindel\"of spaces, and indeed this is done and
answered in \cite{M00}. Indeed \cite{M02} and \cite{M00} combined
cover this question so tightly that the only gap left is that the
example in \cite{M02} is pseudocompact (i.e.~every real function
is bounded) while the example in \cite{M00} is not so, which is exactly
question 1 in \cite{M00}. The simplest example of the spaces discussed
here, $\{0,1\}^{\kappa}\setminus\{\overrightarrow{pt}\}$ fulfills
these conditions --- Tychonoff, discretely star-Lindel\"of, pseudocompact
and with extent $\kappa$. 

The paper is organized as follows: section 2 will be devoted to star-Lindel\"of
spaces, for which the proof is somewhat simpler and also, it seems,
more flexible. Section 3 will be devoted to discretely star-Lindel\"of
spaces. In section 4 we shall give another proof that a $T_{4}$ star-Lindel\"of
space has extent $\leq2^{\aleph_{0}}$ --- although a stronger claim
was proved in \cite{M02}, this proof is simpler and sheds more light
on the examples of chapters 2 and 3.

I wish to thank Mikhail Matveev for encouraging me to publish these
results, for reviewing preprints, and for generally being a nice person.

\section{Star-Lindel\"of spaces}

\begin{thm}
\label{thm_sl}Let $B$ be a compact space and $0\in B$ a point satisfying
that $B\setminus\{0\}$ is Lindel\"of. Then for every cardinality
$\kappa$, the space \[
X=B^{\kappa}\setminus\{\vec{0}\}\]
 is star-Lindel\"of, where $\vec{0}$ is the point all whose coordinates
are $0$.
\end{thm}
\begin{proof}
Let $\Us$ be an open cover of $X$. Without loss of generality, we
may assume $\Us$ is built of basic open sets, i.e.~sets of the form
\[
\{ x\in X\,:\, x(\alpha)\in O_{\alpha},\:\alpha\in E\}\]
where $E\subset\kappa$ is a finite set and the $O_{\alpha}$'s are
open in $B$. We may also assume that for every $U\in\Us$, \begin{equation}
\exists\alpha\in E,\:0\notin O_{\alpha}\label{U_nontrivial}\end{equation}
since otherwise we can simply find a finite subcover of $\Us$. The
first step is to color $\kappa$ by $2^{\aleph_{0}}$ colors, i.e.~to
construct a mapping $c\,:\,\kappa\rightarrow\{0,1\}^{\NN}$ using
the following inductive process: let $C$ be the set of already colored
elements, and let $\alpha\in\kappa$ be $\min\kappa\setminus C$.
The set \[
\{ x\in X\,:\, x(\alpha)\neq0\}\]
 is a product of a Lindel\"of space and a compact space, so it is
Lindel\"of. We take a countable subcover $\{ U_{i}^{\alpha}\}_{i=1}^{\infty}\subset\Us$.
If \[
U_{i}^{\alpha}=\{ x\in X\,:\, x(\beta)\in O_{i,\beta}^{\alpha},\:\beta\in E_{i}^{\alpha}\}\]
 then we define the total index set \[
E^{\alpha}:=\bigcup_{i}E_{i}^{\alpha}\quad.\]
We repeat this process for every uncolored element of $D_{1}^{\alpha}:=E^{\alpha}$
defining\[
D_{n}^{\alpha}:=\bigcup\{ E^{\beta}\,:\,\beta\in D_{n-1}^{\alpha}\setminus C\}\]
and\[
D_{\infty}^{\alpha}:=\bigcup_{n=1}^{\infty}D_{n}^{\alpha}\quad.\]
Clearly, $|D_{\infty}^{\alpha}|\leq\aleph_{0}$. Since $c(D_{\infty}^{\alpha}\cap C)$
is countable we may color $D_{\infty}^{\alpha}\setminus C$ by different
colors taken from $\{0,1\}^{\NN}\setminus c(D_{\infty}^{\alpha}\cap C)$.
This defines the coloring process.

Next, for each $\beta\in D_{\infty}^{\alpha}\setminus C$ we take,
for every $i\in\NN$ an arbitrary element $u_{\beta,i}\in U_{i}^{\beta}$;
and for every $\gamma\in D_{\infty}^{\alpha}$ we define a subset
of $B$, \[
V^{\gamma}:=\{0\}\cup\{ u_{\beta,i}(\gamma)\,:\,\beta\in D_{\infty}^{\alpha}\setminus C,\, i\in\NN\}\quad,\]
{}``the subset of needed values'', which is again countable. Finally,
we want a {}``coloring order'' function $\varphi$, so we define
\[
\varphi(\beta):=\alpha\quad\forall\beta\in D_{\infty}^{\alpha}\setminus C\quad.\]

Recapitulating, we have inductively defined the following objects:
\begin{enumerate}
\item A mapping $c\,:\,\kappa\rightarrow\{0,1\}^{\NN}$.
\item A mapping $\varphi\,:\,\kappa\rightarrow\kappa$.
\item For every $\alpha\in\kappa$ sets $U_{i}^{\alpha}\in\Us$ covering
$\{ x\in X\,:\, x(\alpha)\neq0\}$ and their total index set $E^{\alpha}$.
\item For every $\alpha\in\kappa$ a countable set $V^{\alpha}=\{ v^{\alpha}(1),\, v^{\alpha}(2),\,...\}\subset B$,
$v^{\alpha}(1)=0$.
\end{enumerate}
and these four objects are connected by the facts that \begin{eqnarray}
\beta\in E^{\alpha} & \Rightarrow & \varphi(\beta)\leq\varphi(\alpha)\\
\beta,\gamma\in E^{\alpha}\wedge(\varphi(\beta)=\varphi(\alpha)\vee\varphi(\gamma)=\varphi(\alpha)) & \Rightarrow & c(\beta)\neq c(\gamma)\label{colorful}\\
U_{i}^{\alpha}=\{ x\in X\,:\, x(\beta)\in O_{i,\beta}^{\alpha},\,\beta\in E^{\alpha}\}\wedge\varphi(\beta)=\varphi(\alpha) & \Rightarrow & V^{\beta}\cap O_{i,\beta}^{\alpha}\neq\emptyset\label{Vdense}\end{eqnarray}

The second step is to define the countable set $F$. For every finite
set $I\subset\NN$ and for any function $g\,:\,\{0,1\}^{I}\rightarrow\NN$
we define an element $f_{I,g}$ of $B^{\kappa}$ by\begin{equation}
f_{I,g}(\alpha):=v^{\alpha}(g(\{ c(\alpha)_{i}\}_{i\in I}))\label{def_fIg}\end{equation}
and then \[
F:=\{ f_{I,g}\,:\, f_{I,g}\not\equiv0\}\quad.\]
Clearly, $|F|\leq\aleph_{0}$.

The final step is to show that $\St(F,\Us)=X$. Let therefore $x$
be in $X$, and let $\alpha\in\kappa$ be an element satisfying $x(\alpha)\neq0$
with minimal $\varphi(\alpha)$ (if more than one exists, choose any).
Using the fact that the sets $U_{1}^{\alpha},...$ cover the set $\{ x(\alpha)\neq0\}$
we find some $i$ for which $x\in U_{i}^{\alpha}$. We can represent
$U=U_{i}^{\alpha}$ as \[
U=\{ x\in X\,:\, x(\beta)\in O_{\beta},\:\forall\beta\in E\}\]
where $E\subset E^{\alpha}$ is a finite set.

Now, for every $\beta\in E$ we have some $j_{\beta}\in\NN$ satisfying
\begin{equation}
v^{\beta}(j_{\beta})\in O_{\beta}\label{def_jbeta}\end{equation}
 because $\varphi(\beta)<\varphi(\alpha)$ implies $x(\beta)=0$ so
$0\in O_{\beta}$ and we can choose $j_{\beta}=1$, while for $\varphi(\beta)=\varphi(\alpha)$
we use (\ref{Vdense}). We now write\[
E=C\cup D\]
 with\begin{eqnarray*}
C & := & \{\beta\in E\,:\,\varphi(\beta)<\varphi(\alpha)\}\\
D & := & \{\beta\in E\,:\,\varphi(\beta)=\varphi(\alpha)\}\quad.\end{eqnarray*}
For every $\beta,\gamma\in E$ not both in $C$ we know (\ref{colorful})
that $c(\beta)\neq c(\gamma)$ and we choose an index $i=i_{\beta,\gamma}$
such that\[
c(\beta)_{i}\neq c(\gamma)_{i}\quad.\]
We let $I=\{ i_{\beta,\gamma}\}$ and define a function $g:\{0,1\}^{I}\rightarrow\NN$
using\[
g(\{\epsilon_{i}\})=\left\{ \begin{array}{ll}
j_{\beta} & \exists\beta\in E,\:\epsilon_{i}\equiv c(\beta)_{i}\\
1 & \forall\beta\in E,\:\epsilon_{i}\not\equiv c(\beta)_{i}\end{array}\right.\quad.\]
This is a good definition because the sets $\{ c(\beta)_{i}\}$ can
be identical for $\beta\neq\gamma$ only if $\beta,\gamma\in C$ but
in this case $j_{\beta}=j_{\gamma}=1$. This immediately implies (remember
(\ref{def_fIg}, \ref{def_jbeta})) that \[
f_{I,g}(\beta)\in O_{\beta}\quad\forall\beta\in E\quad.\]
 This implies that $f_{I,g}\not\equiv0$, using (\ref{U_nontrivial}),
thus $f_{I,g}\in F$. Of course, this also gives $f\in U$ so $U\subset\St(F,\Us)$
and since $x\in U$ the theorem is proved.
\end{proof}
If the base space $B$ is $T_{3\frac{1}{2}}$ then the space $X$
defined above is also $T_{3\frac{1}{2}}$. If $B$ contains one point
separated from $0$ (denote this point by $1$) then $e(X)\geq\kappa$
since the set\[
\{ x\in X\,:\,\exists\alpha,\, x(\alpha)=1,\, x(\beta)=0\quad\forall\beta\neq\alpha\}\]
is a closed discrete set. If, say, the topology of $B$ has a base
with size $\leq\kappa$ then $e(X)=\kappa$ since the topology of
$X$ will have a $\kappa$ sized base. Thus, for example, the spaces
$\{0,1\}^{\kappa}\setminus\{\overrightarrow{pt}\}$ is a $T_{3\frac{1}{2}}$
star-Lindel\"of space with $e(X)=\kappa$.

It is interesting to note that removing a little more from $\{0,1\}^{\kappa}$
will destroy this construction. For example, \begin{equation}
Y^{\kappa}:=\{0,1\}^{\kappa}\setminus\{ x\,:\, x(\alpha)=0\:\mathrm{except}\:\mathrm{for}\:\mathrm{at}\:\mathrm{most}\:\mathrm{one}\:\alpha\}\label{def_Yk}\end{equation}
is not star-Lindel\"of for $\kappa>2^{\aleph_{0}}$. We shall present
the proof of this fact in section 4.

Theorem \ref{thm_sl} can be generalized to arbitrary cardinals as
follows:

\begin{thm}
\label{thm:lambda}If $\lambda\leq\tau$, $B$ is a compact space
and $0\in B$ a point such that $B\setminus\{0\}$ is $\lambda$-Lindel\"of then
the star-Lindel\"of number of $B^{\kappa}\setminus\{\vec{0}\}$ is
$\leq\tau$ for any cardinality $\kappa$.
\end{thm}
Where the definitions of $\lambda$-Lindel\"of and the star-Lindel\"of
number are the natural ones: if every cover of $X$ has a subcover
of cardinality $\leq\lambda$ then $X$ is $\lambda$-Lindel\"of,
and the \label{def_starlinnum}star-Lindel\"of number $\tau$ is
the minimal cardinality such that for every cover $\Us$ one has a
set $F$, $|F|\leq\tau$ with $\St(F,\Us)=X$. We shall omit the proof
of this theorem and contend ourselves with the following remark: the
{}``coloring'' step only requires $\lambda<2^{\tau}$ --- the stronger
condition $\lambda\leq\tau$ is necessary for the definition of the
sets $V^{\beta}$. 

Finally, I wish to note that theorem \ref{thm_sl} generalizes without
any chan\-ge to the case of different product terms:

\begin{thm}
Let $B^{\alpha}$, $\alpha\in\kappa$ be compact spaces and let $0^{\alpha}\in B^{\alpha}$
be points satisfying that $B^{\alpha}\setminus\{0^{\alpha}\}$ is
Lindel\"of. Then the space \[
\left(\prod_{\alpha\in\kappa}B^{\alpha}\right)\setminus\{\vec{0}\}\]
is star-Lindel\"of.
\end{thm}
Theorem \ref{thm:lambda} may be generalized in the same manner.

\section{Discretely star-Lindel\"of spaces}

If $B$ is a trivial space, that is only $B$ and $\emptyset$ are
open, say with 2 points, then $B^{\kappa}\setminus\{\vec{0}\}$ has
no non-empty closed discrete subsets. To avoid such issues we shall
restrict our attention to $T_{1}$ spaces.

\begin{thm}
Let $B$ be a $T_{1}$ space, and let $0\in B$ be a point such that
$B\setminus\{0\}$ is a compact space. Then $X:=B^{\kappa}\setminus\{\vec{0}\}$
is discretely star-Lindel\"of.
\end{thm}
Be forewarned that this proof is even messier than the proof of the
previous theorem! The basic idea is the same but there is no {}``coloring''
step and we construct $F$ directly.

\begin{proof}
As before, we shall require from the cover $\Us$ that it is done
with basic open sets, and that it is not trivial (i.e.~no finite
subcover). We shall construct inductively the following objects:
\begin{enumerate}
\item A set $C^{\alpha}\subset\kappa$, {}``the set of fully defined indexes''.
\item A sequence of sets $P_{m}^{\alpha}\subset\kappa\setminus C^{\alpha}$
with $\bigcup_{m}P_{m}^{\alpha}$ finite, {}``the sets of partially
defined indexes'' (in essence at most two will be non-empty at each
step).
\item \label{req_fn_zero}Functions $f_{n}^{\alpha}\,:\, C^{\alpha}\rightarrow B$
satisfying $\forall\beta\in C^{\alpha}$ that $f_{n}^{\alpha}(\beta)=0$
for almost all $n$. The $f_{n}^{\alpha}$'s extend each other i.e.~$\alpha_{1}<\alpha_{2}$
implies $C^{\alpha_{1}}\subset C^{\alpha_{2}}$ and $\left.f_{n}^{\alpha_{2}}\right|_{C^{\alpha_{1}}}\equiv f_{n}^{\alpha_{1}}$.
\item For every $\beta\in\bigcup_{m}P_{m}^{\alpha}$ a finite set $I_{\beta}^{\alpha}\subset\NN$
(we assume as a matter of notation that $\beta\notin\cup P_{m}^{\alpha}$
implies $I_{\beta}^{\alpha}=\emptyset$) and scalars $\overline{f_{n}^{\alpha}}(\beta)\in B$
for every $n\in I_{\beta}^{\alpha}$.
\item A one-to-one function $\varphi^{\alpha}\,:\, C^{\alpha}\rightarrow[1,\alpha]$
such that $\varphi(\beta)$ is always a non-limit ordinal, {}``the
defining order function'', with the $\varphi^{\alpha}$'s also extending
each other.
\end{enumerate}
Assume $\alpha$ is non-limit, i.e.~$\alpha^{-}$ satisfies that
$\left(\alpha^{-}\right)^{+}=\alpha$. First we select an $m^{\alpha}\in\NN$
and a $\mu^{\alpha}\in\kappa\setminus C^{\alpha^{-}}$ as follows:
if for some $m$ we have $P_{m}^{\alpha^{-}}\neq\emptyset$ then we
take $m^{\alpha}$ to be the minimal such $m$ and $\mu^{\alpha}:=\min P_{m}^{\alpha^{-}}$.
Otherwise, we take $m^{\alpha}:=0$ and $\mu^{\alpha}:=\min\kappa\setminus C^{\alpha^{-}}$.

Now we examine the set \[
X_{\mu^{\alpha}}:=\{ x\in X\,:\, x(\mu^{\alpha})\neq0\}\]
which is compact, so we take a finite subcover $U_{1}^{\alpha},...,U_{k^{\alpha}}^{\alpha}$
and let $E^{\alpha}$ be the total index set. For every $\beta\in E^{\alpha}\cap C^{\alpha^{-}}$
we define\begin{eqnarray*}
J_{\beta}^{\alpha} & := & \{ n\in\NN\,:\, f_{n}^{\alpha^{-}}(\beta)\neq0\}\\
J^{\alpha} & := & \bigcup J_{\beta}^{\alpha}\end{eqnarray*}
which is a finite set. Finally, for every $i\leq k^{\alpha}$ we take
some $u_{i}^{\alpha}\in U_{i}^{\alpha}$ and a distinct number \begin{equation}
n_{i}^{\alpha}\in\NN\setminus(J^{\alpha}\cup\bigcup_{\beta\in\cup P_{m}^{\alpha^{-}}}I_{\beta}^{\alpha^{-}})\quad.\label{def_ni}\end{equation}
We are now ready to proceed with the induction. We define\begin{align}
C^{\alpha} & :=C^{\alpha^{-}}\cup\{\mu^{\alpha}\}\nonumber \\
f_{n}^{\alpha}(\beta) & :=\left\{ \begin{array}{ll}
f_{n}^{\alpha^{-}}(\beta) & \beta\neq\mu^{\alpha}\\
\overline{f_{n}^{\alpha^{-}}}(\beta) & \beta=\mu^{\alpha},\, n\in I_{\mu^{\alpha}}^{\alpha^{-}}\\
u_{i}^{\alpha}(\beta) & \beta=\mu^{\alpha},\, n=n_{i}^{\alpha}\\
0 & \mathrm{otherwise}\end{array}\right.\label{def_fn}\\
P_{m}^{\alpha} & :=\left\{ \begin{array}{ll}
\left(P_{m}^{\alpha^{-}}\setminus\{\mu^{\alpha}\}\right)\cup\left(E^{\alpha}\setminus C^{\alpha}\right) & m=m^{\alpha}+1\\
P_{m}^{\alpha^{-}}\setminus\{\mu^{\alpha}\} & \mathrm{otherwise}\end{array}\right.\label{def_Pm}\\
I_{\beta}^{\alpha} & :=\left\{ \begin{array}{ll}
I_{\beta}^{\alpha^{-}}\cup\{ n_{i}^{\alpha}\}_{i=1}^{k^{\alpha}} & \beta\in E^{\alpha}\setminus C^{\alpha}\\
I_{\beta}^{\alpha^{-}} & \mathrm{otherwise}\end{array}\right.\label{def_Ibeta}\\
\overline{f_{n}^{\alpha}}(\beta) & :=\left\{ \begin{array}{ll}
u_{i}^{\alpha}(\beta) & \beta\in E^{\alpha}\setminus C^{\alpha},\, n=n_{i}^{\alpha}\\
\overline{f_{n}^{\alpha^{-}}}(\beta) & \mathrm{otherwise}\end{array}\right.\label{def_fnbar}\\
\varphi^{\alpha}(\beta) & :=\left\{ \begin{array}{ll}
\alpha & \beta=\mu^{\alpha}\\
\varphi^{\alpha^{-}}(\beta) & \mathrm{otherwise}\end{array}\right.\nonumber \end{align}
which clearly fulfills all finiteness requirements.

For $\alpha$ a limit ordinal we define more simply\begin{eqnarray}
C^{\alpha} & := & \bigcup_{\beta<\alpha}C^{\beta}\nonumber \\
P_{m}^{\alpha} & := & \emptyset\quad\forall m\label{Blimit_empty}\\
f_{n}^{\alpha} & := & \bigcup_{\beta<\alpha}f_{n}^{\beta}\nonumber \\
\varphi^{\alpha} & := & \bigcup_{\beta<\alpha}\varphi^{\beta}\quad.\nonumber \end{eqnarray}
This finishes the description of the induction, and we must now show
that it actually creates relevant objects. We start off with something
light.
\begin{claim*}
$C^{\kappa}=\kappa$
\end{claim*}
We build inductively a $\psi:\kappa\rightarrow\kappa$ such that $[1,\alpha)\subset C^{\psi(\alpha)}$
and $|\psi(\alpha)|=|\alpha|$ for every infinite $\alpha$. For a
non-limit ordinal $\alpha$ we take \[
\psi(\alpha):=\psi(\alpha^{-})+\omega+1\]
and the induction hypothesis is fulfilled due to (\ref{Blimit_empty})
and to the definition of $m^{\alpha}$, while for a limit $\alpha$
we take\[
\psi(\alpha):=\bigcup_{\beta<\alpha}\psi(\beta)\]
which clearly finishes the claim. 

This shows that $f_{n}^{\kappa}$ and $\varphi^{\kappa}$ are indeed
well defined functions on $\kappa$, so define\begin{eqnarray*}
f_{n} & := & f_{n}^{\kappa}\\
\varphi & := & \varphi^{\kappa}\\
F & := & \{ f_{n}\,:\, f_{n}\not\equiv0\}\end{eqnarray*}
The fact that $F$ is closed and discrete follows immediately from
requirement \ref{req_fn_zero}, i.e.~from \[
f_{n}(\beta)=0\:\mathrm{for}\:\mathrm{almost}\:\mathrm{all}\: n\]
(here we used the fact that $B$ is $T_{1}$) so we are now left with
the chore of showing that $\St(F,\Us)=X$. Let $x$ be in $X$ and
let \begin{eqnarray*}
\alpha & := & \min\{\varphi(\beta)\,:\, x(\beta)\neq0\}\\
\mu^{\alpha} & = & \varphi^{-1}(\alpha)\quad.\end{eqnarray*}
This definition of $\alpha$ implies that it is a non-limit ordinal,
and that $x(\beta)=0$ for all $\beta\in C^{\alpha^{-}}$. We pick
some $i(x)$ such that $x\in U_{i(x)}^{\alpha}$, and again write\[
U_{i(x)}^{\alpha}=:\{ u\in X\,:\, u(\beta)\in O_{\beta},\:\forall\beta\in E^{\alpha}\}\quad.\]
We examine the corresponding $n_{i(x)}^{\alpha}$. Denote $n=n_{i(x)}^{\alpha}$
and $f=f_{n}$.
\begin{enumerate}
\item For $\beta\in E^{\alpha}\cap C^{\alpha^{-}}$ we have $n\notin J_{\beta}^{\alpha}$
(remember (\ref{def_ni})) so $f(\beta)=0$ and since $x(\beta)=0$
we must have $f(\beta)\in O_{\beta}$.
\item For $\beta=\mu^{\alpha}$ we have $f(\beta)=u_{i(x)}^{\alpha}(\beta)$
((\ref{def_fn}) clause 3) and since $u_{i(x)}^{\alpha}\in U_{i(x)}^{\alpha}$
we again get $f(\beta)\in O_{\beta}$.
\item For $\beta\in E^{\alpha}\setminus C^{\alpha}$ we have ((\ref{def_Pm}),
(\ref{def_Ibeta}) and (\ref{def_fnbar})) that $\beta\in P_{m^{\alpha}+1}^{\alpha}$,
$n\in I_{\beta}^{\alpha}$ and \[
\overline{f_{n}^{\alpha}}(\beta)=u_{i(x)}^{\alpha}(\beta)\in O_{\beta}\quad.\]
Now the definition of $m^{\alpha}$ gives that $P_{m}^{\alpha}=\emptyset$
for all $m<m^{\alpha}$, and the definition of $\mu^{\alpha}$ gives
that every step from $\alpha$ to $\alpha^{+}$ removes one element
from $P_{m^{\alpha}}$ so for $\alpha':=\alpha+\left|P_{m^{\alpha}}^{\alpha}\right|$
we have $P_{m^{\alpha}}^{\alpha'}=\emptyset$ as well and for $\alpha'':=\alpha'+\left|P_{m^{\alpha}+1}^{\alpha'}\right|$
we have $P_{m^{\alpha}+1}^{\alpha''}=\emptyset$ and in particular
$\beta\in C^{\alpha''}$ so \[
f_{n}^{\alpha''}(\beta)=\overline{f_{n}^{\alpha}}(\beta)\]
((\ref{def_fn}) clause 2 for some appropriate $\alpha'\leq\overline{\alpha}\leq\alpha''$)
so again $f(\beta)\in O_{\beta}$.
\end{enumerate}
These 3 give us that $f(\beta)\in O_{\beta}$ for all $\beta\in E^{\alpha}$
so $f\in U_{i(x)}^{\alpha}$. As before, we can conclude that $f\not\equiv0$
since otherwise $\Us$ would have a finite subcover and therefore
$f\in F$ which means that $x\in U_{i(x)}^{\alpha}\subset\St(F,\Us)$
and the theorem is proved.
\end{proof}
As explained in section 2, the space $X:=\{0,1\}^{\kappa}\setminus\{ pt\}$
is Tychonoff and $e(X)=\kappa$. It is well known that this space
is pseudocompact therefore it is also an answer to question 1 from
\cite{M00}. I provide a proof for the convenience of the reader.

\begin{thm}
$X=\{0,1\}^{\kappa}\setminus\{ pt\}$ is pseudocompact for every $\kappa>\aleph_{0}$.
\end{thm}
\begin{proof}
Let $f\,:\, X\rightarrow\RR$ be a continuous unbounded function.
We take $U_{n}$ to be basic open sets inside $f^{-1}(\RR\setminus[-n,n])$.
We write\[
U_{n}=\{ x\in X\,:\, x(\beta)\in O_{\beta},\:\forall\beta\in E_{n}\}\quad.\]
We find a $g_{1}\,:\, F_{1}\rightarrow\{0,1\}$ ($F_{1}:=E_{1}$),
such that infinitely many $U_{n}$'s intersect the set \[
X_{1}:=\{ x\in X\,:\, x|_{F_{1}}\equiv g_{1}\}\]
(this is clearly possible). Continuing inductively, we take $n_{2}$
to be some index satisfying $U_{n_{2}}\cap X_{1}\neq\emptyset$ define
$F_{2}:=F_{1}\cup E_{n_{2}}$ and $g_{2}\,:\, F_{2}\rightarrow\{0,1\}$
an extension of $g_{1}$ satisfying that infinitely many $U_{n}$'s
intersect $X_{2}$. Continuing this process $\aleph_{0}$ steps we
reach a contradiction since any element of $X$ such that \[
x|_{\cup F_{i}}=\cup g_{i}\]
is a point of discontinuity for $f$.
\end{proof}
\begin{rem*}
Another popular variation on the definition is to require that the
set $F$ will be taken in some predefined dense subset $S$. See for
example \cite{B98}, for the properties {}``absolutely star-Lindel\"of'',
in which $F$ is assumed to be countable; and {}``property (a)'',
in which $F$ is assumed to be closed and discrete, but not necessarily
countable. Nothing like that works for our spaces, e.g.~for $\{0,1\}^{\kappa}\setminus\{\vec{0}\}$,
$\kappa>\aleph_{0}$. For absolutely star-Lindel\"of spaces, this
is proved in \cite[theorem 8.2]{B98}. The following simple proof
shows this for both definitions. Take as the dense set \[
S:=\{ x\in X\,:\, x(\alpha)=1\:\mathrm{except}\:\mathrm{on}\:\mathrm{a}\:\mathrm{finite}\:\mathrm{set}\}\quad.\]
For the cover, divide $\kappa$ into $\kappa_{1}\cup\kappa_{2}$,
$|\kappa_{1}|=|\kappa_{2}|=\kappa$ and let $\varphi\,:\,\kappa_{1}\rightarrow\kappa_{2}$
be one-to-one onto. Define \[
D:=\{ x\in X\,:\,\exists\alpha\in\kappa_{1},\, x(\beta)=1\Leftrightarrow\alpha=\beta\}\]
which is closed so we can cover $X$ using sets of the form \[
U_{\alpha}:=\{ x\,:\, x(\alpha)=1,\: x(\varphi(\alpha))=0\},\quad\alpha\in\kappa_{1}\]
and add $X\setminus D$. This is obviously an irreducible cover, that
is it has no proper subcover. To show that $X$ is not {}``absolutely
star-Linel\"of'', take a countable $F\subset S$. Clearly for some
$\alpha\in\kappa_{1}$, $f(\alpha)=f(\varphi(\alpha))=1$ for all
$f\in F$  and therefore $\St(F,\Us)\neq X$.

As for closed and discrete subsets of $S$ ({}``property (a)''),
note just that $S$ contains only finite closed discrete subsets,
so this problem reduces to the previous one.
\end{rem*}

\section{$T_{4}$ spaces}

\begin{thm}
Any $T_{4}$ star-Lindel\"of space satisfies $e(X)\leq2^{\aleph_{0}}$.
\end{thm}
\begin{proof}
Assume to the contrary that $D\subset X$ is a closed discrete subset,
$|D|=\kappa>2^{\aleph_{0}}$. Let $\sigma$ be a one-to-one onto mapping\[
\sigma\,:\, D\rightarrow\{(\alpha,\beta)\in\kappa\times\kappa\,:\,\alpha\neq\beta\}\quad.\]
For every $\alpha\in\kappa$ we define the closed sets\begin{eqnarray*}
A_{\alpha} & := & \sigma^{-1}(\{(\alpha,\beta)\,:\,\beta\neq\alpha\})\\
B_{\alpha} & := & \sigma^{-1}(\{(\beta,\alpha)\,:\,\beta\neq\alpha\})\quad.\end{eqnarray*}
Using normality, we find two disjoint open sets $C_{\alpha}$ and
$D_{\alpha}$ satisfying\begin{eqnarray*}
C_{\alpha}\cap D & = & A_{\alpha}\\
D_{\alpha}\cap D & = & B_{\alpha}\quad.\end{eqnarray*}
We examine the family \[
\Us:=\{ C_{\alpha}\cap D_{\beta}\,:\,\alpha\neq\beta\}\]
adding to it, if necessary, the set $X\setminus D$ to make it a cover.
This is clearly an irreducible cover. Assume that for some countable
$F$, $\St(F,\Us)=X$ . We write $F=\{ f_{1},f_{2},...\}$ and build
a map\[
\tau\,:\,\kappa\rightarrow\{0,1\}^{\NN}\]
 by\[
\tau(\alpha)_{i}=\left\{ \begin{array}{ll}
1 & f_{i}\in D_{\alpha}\\
0 & \mathrm{otherwise}\end{array}\right.\quad.\]
Since $\kappa>2^{\aleph_{0}}$ we have some $\alpha\neq\beta$ for
which $\tau(\alpha)=\tau(\beta)$. However, in this case, \[
f_{i}\notin C_{\alpha}\cap D_{\beta}\quad\forall i\in\NN\]
 since \[
f_{i}\in C_{\alpha}\cap D_{\beta}\Rightarrow\tau(\alpha)_{i}\neq\tau(\beta)_{i}\]
and therefore \begin{align*}
\sigma^{-1}((\alpha,\beta)) & \notin\St(F,\Us)\qedhere\end{align*}

\end{proof}
As discussed in the introduction, this result is not new. However,
we will now use this technique to show, as promised, that the spaces
$Y^{\kappa}$ defined in section 2 (equation (\ref{def_Yk})) are
not star-Lindel\"of for $\kappa>2^{\aleph_{0}}$. We start by dividing
$\kappa$ into $\kappa_{1}\cup\kappa_{2}$ with $|\kappa_{1}|=|\kappa_{2}|=\kappa$,
taking a one-to-one onto map \[
\sigma\,:\,\kappa_{1}\rightarrow\{(\alpha,\beta)\in\kappa_{2}\times\kappa_{2}\,:\,\alpha\neq\beta\}\quad,\]
a closed discrete set \[
D=\{ x\in Y^{\kappa}\,:\, x(\alpha)=1\:\mathrm{only}\:\mathrm{for}\:\mathrm{one}\:\alpha\in\kappa_{1}\:\mathrm{and}\:\mathrm{for}\:\sigma(\alpha)_{1}\}\quad,\]
 open sets\[
U_{\alpha}:=\{ x\in Y^{\kappa}\,:\, x(\alpha)=1,\: x(\sigma(\alpha)_{1})=1,\: x(\sigma(\alpha)_{2})=0\}\quad,\]
and a cover\[
\Us:=\{ U_{\alpha}\}_{\alpha\in\kappa_{1}}\cup\{ Y^{\kappa}\setminus D\}\]
which is clearly irreducible. The proof now continues similarly, with
$C_{\alpha}$ and $D_{\alpha}$ defined for $\alpha\in\kappa_{2}$
by \begin{eqnarray*}
C_{\alpha} & := & \{ x\in Y^{\kappa}\,:\, x(\alpha)=0\}\\
D_{\alpha} & := & \{ x\in Y^{\kappa}\,:\, x(\alpha)=1\}\quad.\end{eqnarray*}


\begin{thebibliography}{BM01}
\bibitem[B98]{B98}Maddalena Bonanzinga, \emph{Star-Lindelof and absolutely star-Lindelof
spaces,} Q \& A in General Topology, 16 (1998), 79-104.
\bibitem[BM01]{BM01}Maddalena Bonanzinga and Mikhail V. Matveev, \emph{Products of Star-Lindel\"of
and Related Spaces}, Houston Journal of Mathematics, 27/1 (2001),
45-57.
\bibitem[M98]{M98}Mikhail V. Matveev, \emph{A survey on star covering properties}, Topological
Atlas, Preprint 330, \texttt{at.yorku.ca/v/a/a/a/19.htm}
\bibitem[M00A]{M02}Mikhail V. Matveev, \emph{How weak is weak extent?}, Topology and
its Applications, 119/2 (2002), 229-232. \texttt{arXiv.org/abs/math.GN/0006198}
\bibitem[M00B]{M00}Mikhail V. Matveev, \emph{On Spaces in Countable Web}, \texttt{\url{arXiv.org/abs/math.GN/0006199}}
\bibitem[YG99]{YY99}Yoshikazu Yasui and Zhi-min Gao, \emph{Spaces in Countable Web}, Houston
Journal of Mathematics, 25 (1999), 327-330.
\end{thebibliography}
\end{document}